\documentclass[11pt]{article}

\usepackage[a4paper,margin=2.5cm]{geometry}

\usepackage{amsfonts}
\usepackage{amssymb}
\usepackage{amsmath}
\usepackage{bm}
\usepackage{graphicx}
\usepackage[colorlinks=true, allcolors=blue]{hyperref}
\usepackage{authblk}
\usepackage{subfigure}
\usepackage{latexsym}
\usepackage{epstopdf}
\usepackage{epsfig}
\usepackage{setspace}
\usepackage{algorithm}
\usepackage{algpseudocode}
\usepackage{booktabs}
\usepackage{multirow}
\usepackage{xurl}

\spacing{1.5}

\graphicspath{{img/}}

\title{Optimal experimental design for parameter estimation \\ in the presence of observation noise}
\author[1]{Jie Qi}
\author[2]{Ruth E. Baker}
\affil[1]{College of Information Science and Technology, Donghua University, Shanghai, China}
\affil[2]{Mathematical Institute, University of Oxford, Oxford, United Kingdom} 
\date{}

\begin{document}
\maketitle


\begin{abstract}
Using mathematical models to assist in the interpretation of experiments is becoming increasingly important in research across applied mathematics, and in particular in biology and ecology. In this context, accurate parameter estimation is crucial; model parameters are used to both quantify observed behaviour, characterise behaviours that cannot be directly measured and make quantitative predictions. The extent to which parameter estimates are constrained by the quality and quantity of available data is known as parameter identifiability, and it is widely understood that for many dynamical models the uncertainty in parameter estimates can vary over orders of magnitude as the time points at which data are collected are varied. Here, we use both local sensitivity measures derived from the Fisher Information Matrix and global measures derived from Sobol’ indices to explore how parameter uncertainty changes as the number of measurements, and their placement in time, are varied. We use these measures within an optimisation algorithm to determine the observation times that give rise to the lowest uncertainty in parameter estimates. Applying our framework to models in which the observation noise is both correlated and uncorrelated demonstrates that correlations in observation noise can significantly impact the optimal time points for observing a system, and highlights that proper consideration of observation noise should be a crucial part of the experimental design process. 
\end{abstract}

\pagebreak


\section{Introduction}

Mathematical modelling serves as a crucial technique in the interpretation of experiments and offers conceptual insights into the mechanisms underlying complex systems across various fields, particularly in biology and ecology. In this context, accurate parameter estimation is essential; model parameters are now routinely employed to quantify observed behaviours, characterize behaviours that cannot be directly measured, and make quantitative predictions. The extent to which parameter estimates are constrained by the quality and quantity of available data is termed parameter identifiability. For many dynamical models, the uncertainty in parameter estimates can vary over orders of magnitude as the time points at which data are collected are varied. Consequently, optimising data collection strategies to improve parameter identifiability is a crucial aspect of experimental design, particularly for biological systems that frequently exhibit complex and non-linear behaviours. 

Optimal experimental design methodologies provide a means to optimise experimental protocols to e.g., maximise information acquisition or minimise uncertainty in parameter estimates under resource constraints. The choice of objective function in optimal experimental design is crucial as these set specific optimality criteria, and sensitivity measures are commonly employed. In particular, local sensitivity measures incorporated in the Fisher information matrix~\cite{gutenkunst3universally,chis2016relationship,telen2014robustifying,zhang2018optimal} are often used because the inverse of the Fisher information matrix provides an estimate of the lower bound of the covariance matrix of the estimated parameters, based on the Cramér-Rao inequality~\cite{walter1997identification}. Local sensitivity measures assume a linear relationship between parameter perturbations and model responses locally in the region of the specified parameter values. A drawback is therefore that local sensitivity measures can lead to inefficient experimental designs if the input parameters are not close to their ``true'' values. To overcome the limitations of local sensitivity measures, global sensitivity analyses have been increasingly adopted in optimal experimental design~\cite{schenkendorf2018impact}. These measures consider nonlinear effects, non-monotonic behaviour and dependencies between multiple parameters~\cite{iooss2015review}, making them more suited to the complexities of biological systems. In addition, global sensitivities take parameter ranges as inputs, enhancing the robustness of the experimental design process~\cite{pozzi2020global,chu2013necessary}. 

A key unexplored question in the literature is the extent to which changes in the characteristics of the measurement / observation noise process impact the output of optimal experimental design protocols. As such, the goal of this paper is to apply optimal experimental design methodologies to a widely used model in mathematical biology to determine optimal observation schemes, specifically focussing on the time points at which data are collected and how different types of  noise influence optimal data collection strategies. We will explore the ability of both local and global measures, used within experimental design approaches, to reduce the uncertainty in parameter estimates in the face of observation noise.

The model we consider is an ordinary differential equation model, as is commonly used to describe dynamical systems in biology, and the observations are assumed to consist of the solution of this  model plus an ``observation noise''. In this work, our assumption is that that the observation noise accounts for factors not explicitly included in the model, for example, stochasticity in the underlying processes or measurement errors. Most modelling studies in mathematical biology assume that the observation noise is independent Gaussian distributed. This assumption simplifies much statistical analysis, but it does not hold universally across all applications~\cite{lambert2023autocorrelated,fu2022quantifying}. Although there has been some research into the impact of autocorrelated observation noise processes on parameter estimation~\cite{lambert2023autocorrelated,kuhlmann2001importance,simoen2013prediction}, the extent to which it affects optimal experimental design remains under-explored. As such, the overarching aim of this work is to understand the extent to which changes in the observation noise process impact our ability to estimate model parameters and the ``optimal'' experimental set-up that minimises the uncertainty in parameter estimates. 

In particular, in many applications there may be correlations between adjacent observation error terms in time series models; these may be caused by, for example, measurement equipment biases or model misspecification~\cite{lambert2023autocorrelated,fu2022quantifying}. Autocorrelated noise can be represented both through discrete-time models, such as autoregressive, moving average, autoregressive moving average and autoregressive integrated moving average processes~\cite{durbin2012time}, and continuous-time processes, such as fractional Brownian motion~\cite{biagini2008stochastic}, Cox-Ingersoll-Ross processes~\cite{dereich2012euler} and Ornstein-Uhlenbeck (OU)~\cite{maller2009ornstein} processes. The OU process, a continuous-time analogue of the discrete-time autoregressive process of order 1, is particularly suited for analysing the impact of autocorrelated noise on parameter estimation because the extent of the autocorrelation between time points depends solely on the magnitude of the time between them. Therefore, in this study, we use the OU process to model autocorrelated observation noise. 
 
This paper explores optimal experimental design within the context of both uncorrelated (independent, identically distributed (IID)) and autocorrelated (OU distributed) observation noise, specifically focussing on the logistic model, which is ubiquitously applied in the modelling of biological systems. We suggest an approach to optimal experimental design to minimise parameter uncertainty, utilising both local sensitivity measures derived from the Fisher information matrix and global sensitivity measures derived from Sobol’ indices as objective functions, and we use our approach to show that changes in the structure of the observation error model can significantly impact the optimal experimental design.

The paper is structured as follows. Section~\ref{section_methods} introduces the relevant methodology, including the logistic model and the information measures and noise processes to be used in this work. It also describes formulation of the optimal experimental design problem and an overview of the profile-likelihood approach for estimating parameters and confidence intervals. The main results, including the impact of the specific noise process on parameter estimates and the optimal observation schemes, are given in Section~\ref{section_results}. The paper concludes with a brief discussion in Section~\ref{section_discussion}.
    

\section{Methods}
\label{section_methods}

In this section, we outline the mathematical model and relevant statistical techniques used throughout the paper. All relevant code to implement the techniques and reproduce the figures in this work is provided on Github at \url{https://github.com/Jane0917728/Optimal-Experimental-Design-and-Parameter-Estimation-with-Autocorrelated-Observation-Noise}.


\subsection{Mathematical model} 
\label{section_model}

We consider the logistic model for population growth,
\begin{equation}
    \frac{\text{d}C(t)}{\text{d}t} = rC(t)\left(1 - \frac{C(t)}{K}\right),
\end{equation}
where $C(t)$ is the population at time $t$, $r>0$ is the growth rate, $K>0$ is the carrying capacity, with $K=\lim_{t\to\infty}C(t)$, and $C(0)=C_0>$ is the initial population size. For $C_0\ll{K}$, population growth is characterised by an initial exponential growth phase in which resources are assumed abundant, a deceleration phase where competition for resources becomes important, and a steady state phase in which the population stabilises at the carrying capacity. The logistic model can be solved explicitly to give
\begin{equation}
    C(t) = \frac{C_0 K}{(K-C_0)e^{-rt} +C_0}.
    \label{solution_C}
\end{equation}
In the model, three parameters $\bm{\theta}=(\theta_1,\theta_2,\theta_3)=(r,K,C_0)$ are to be estimated, and we will often write $C(t;\bm{\theta})$ to highlight the dependence of the model output on the parameters. The ``true'' parameters are denoted by $\bm{\theta}^*=(r^*,K^*,C_0^*)$. Throughout this work, we set $r=0.2$, $K=50$ and $C_0=4.5$ and we consider the model up to a final time $t_\text{final}=80$, which is sufficient for the population to be very close to carrying capacity.

\pagebreak

We assume that it is possible to observe the state of the system at $n_s$ (strictly positive) discrete \emph{observation times}, $t_1,\ldots,t_{n_s}$, and that observations are of the form
\begin{equation}
\label{model_noise}
    \bm{Y}=\bm{C}(\bm{\theta}^*)+\bm{\epsilon},
\end{equation}
where 
\begin{align}
\label{def_vecC}
    \bm{Y}&=[Y(t_1),\ldots,Y(t_{n_s})],\\
    \bm{C}(\bm{\theta}^*)&=[C(t_1;\bm{\theta}^*),\ldots,C(t_{n_s};\bm{\theta}^*)],\\
    \bm{\epsilon}&=[\epsilon(t_1),\ldots,\epsilon(t_{n_s})],\label{def_vec_epsilon}
\end{align}
represent the vector forms of the observations, models outputs and measurement noise at discrete observation times $t_1,\ldots,t_{n_s}$, respectively. 


\subsection{Measures of information} 
\label{section_information}

To quantify the amount of information about the parameters contained in a dataset, we will consider both the Fisher information matrix, which provides a local measure of uncertainty, and Sobol' indices, which provide a global measure of uncertainty. 


\paragraph{Fisher information matrix.} 

The Fisher information matrix can be defined using the expectation of the Hessian of the log-likelihood function as~\cite{gutenkunst3universally,chis2016relationship}
\begin{equation}
\label{definition_FIM}
    \mathcal{F}=\{F_{ij}\}=\left\{-\mathbb{E}\left[\frac{\partial^2 L(\bm \theta|\bm{Y})}{\partial \log(\theta_i) \partial \log(\theta_j)}\right]\right\},
\end{equation}
where $\mathbb{E}[\cdot]$ refers to the expectation operator and $L(\bm{\theta}|\bm{Y})$ is the log-likelihood of parameters $\bm{\theta}$ given data $\bm{Y}$. The inverse of the Fisher information matrix offers a lower bound on the covariance matrix for unbiased estimators via the Cramér-Rao bound~\cite{walter1997identification}. As such, maximising an objective function based on the Fisher information matrix is equivalent to minimising the uncertainty of the estimated parameters governed by the covariance matrix.


\paragraph{Sobol' indices.} 

Sobol' indices are based on the idea of variance decomposition and evaluate the contribution of each parameter's variance to the total variance of the measurements~\cite{schenkendorf2018impact,borgonovo2016sensitivity,saltelli1995sensitivity}. They are particularly suitable for measuring the sensitivity of complex models where multiple parameters influence the measurements, and the relationships between them are non-linear or involve interactions among parameters.  

\pagebreak

Given a model of the form $M=M(t;\bm{\theta})$, where $\bm{\theta}=(\theta_1,\ldots,\theta_p)$ (recall that here that $p=3$ with $\theta_1=r$, $\theta_2=K$ and $\theta_3=C_0$), the total-effect Sobol' index is defined, at time $t$, as
\begin{equation}
\label{Sobol_total}    
S_i(t)=\frac{\mathbb{E}_{\bm{\theta}_{\sim i}}[\mathbb{V}_{\theta_{i}}[M(t;\bm{\theta})|\bm{\theta}_{\sim i}]]}{\mathbb{V}[M(t;\bm{\theta})]}, 
\quad\text{for~} i=1,\ldots,p,
\end{equation}
where $\mathbb{V}$ is the total variance, $\mathbb{V}_{\theta_i}$ denotes the variance over all possible values of $\theta_i$, and  the inner expectation operator $\mathbb{E}_{\bm{\theta}_{\sim i}}[M(t;\bm{\theta})|\theta_i]$ indicates, given $\theta_i$ fixed, the mean of $M(t;\bm{\theta})$ over all possible values of $\bm{\theta}_{\sim i}$, the set of the parameters excluding $\theta_i$. $S_i(t)$ measures the total effect of $\theta_i$ on the model output $M(t;\bm{\theta})$, including both its direct influence and its interactions with other parameters, and it is normalised such that $S_i(t)\in[0,1]$. A high value of $S_i(t)$ indicates that the $i^\text{th}$ parameter has a strong effect on the output. 


\subsection{Independent and identically distributed observation noise} 
\label{section_IID}

For IID Gaussian observation noise we have $\bm{\epsilon}\sim \mathcal{N} (\bm{0},\sigma^2_\text{IID}\bm{I})$, where $\sigma^2_\text{IID}>0$ is the noise variance, so that $\bm{Y}\sim \mathcal{N} (\bm{C}(\bm{\theta}^*),\sigma^2_\text{IID}\bm{I})$ and the log-likelihood function is given by 
\begin{equation}
\label{Log_likelihood}
    L(\bm{\theta}|\bm{Y})=-\frac{n_s}{2}\log(2\pi\sigma^2_\text{IID})-\frac{1}{2\sigma^2_\text{IID}}(\bm{Y}-\bm{C}(\bm{\theta}))(\bm{Y}-\bm{C}(\bm{\theta}))^{\mathrm{T}}.
\end{equation}
Applying maximum likelihood estimation to estimate the noise variance gives
\begin{equation}\label{eq_sigma_estimate}
    \hat\sigma^2_\text{IID} =\frac{1}{n_s}(\bm{Y}-\bm{C}(\bm \theta)) (\bm{Y}-\bm{C}(\bm \theta))^{\mathrm T}=
    \frac{1}{n_s}\sum_{s=1}^{n_s}(Y(t_s)-C(t_s;\bm{\theta}))^2.
\end{equation}
Throughout this work we will fix the estimate for $\sigma^2_\text{IID}$ to be the maximum likelihood estimate (MLE) $\hat{\sigma}^2_\text{IID}$ given in Equation \eqref{eq_sigma_estimate}, and focus on optimising the measurement times to best estimate the model parameters $\bm{\theta}$.

To calculate the Fisher information matrix, $\mathcal{F}$, we first calculate 
\begin{equation}
     \frac{\partial^2  L(\bm{\theta}|\bm{Y})}{\partial \log(\theta_i)\partial \log(\theta_j)} = -\frac{1}{ \sigma^2_\text{IID}}\theta_i\theta_j \frac{\partial \bm{C}}{\partial \theta_i}   \left(\frac{\partial \bm C}{\partial \theta_j}\right)^{\mathrm{T}}+\frac{1}{\sigma^2_\text{IID}}\theta_i\theta_j(\bm{Y}-\bm{C}(\bm{\theta})) \left(\frac{\partial^2 \bm C}{\partial \theta_i\partial \theta_j}\right)^{\mathrm{T}}.
\end{equation}
Exploiting the fact that $\mathbb{E}[\bm{Y}-\bm{C}(\bm \theta)]=0$, the elements of the Fisher information matrix can be written
\begin{equation}
\label{definition_FIM_2}
     F_{ij}=\frac{1}{\hat{\sigma}^2_\text{IID}}\theta_i\theta_j\frac{\partial \bm{C}}{\partial \theta_i}\left(\frac{\partial \bm{C}}{\partial \theta_j}\right)^{\mathrm{T}}.
\end{equation}

\pagebreak

We construct a global information matrix, $\mathcal{G}$, for IID Gaussian observation noise based on the total effect Sobol' indices with entries
\begin{equation}
\label{def_IIDSob}
     G_{ij}=\bm{S}_{i}\bm{S}_{j}^\mathrm{T}=\sum_{s=1}^{n_s} S_{i}(t_s) S_{j}(t_s),
\end{equation}
where $\bm{S}_{i}=[S_i(t_1),\ldots, S_i(t_{n_s})]$ is a row vector composed of the total effect Sobol' indices of parameter $i$ at each observation time $t_s$, $s=1,\ldots,n_s$, as defined in Equation \eqref{Sobol_total}.  Sobol' indices are usually computed using Monte Carlo simulation~\cite{pozzi2020global,saltelli1995sensitivity}. In this paper, we  compute the Sobol' indices directly using the Global Sensitivity Analysis Toolbox for Matlab~\cite{cannavo2012sensitivity}.


\subsection{Autocorrelated observation noise} 
\label{section_OU}

We use an OU process, as defined by the stochastic differential equation
\begin{equation}
\label{SDE-OU}
    \text{d}\epsilon(t)=-\phi\,\epsilon(t)\text{d}t+\sigma_\text{OU}\text{d}W(t),\quad \epsilon(0)=\epsilon_0,
\end{equation}
to model autocorrelated measurement noise, with $\epsilon_0 \sim \mathcal{N}(0,\sigma^2_\text{OU} /(2\phi))$. In Equation~\eqref{SDE-OU}, $\phi>0$ is the mean-reversion rate and $\sigma_\text{OU}>0$, often known as the volatility coefficient, controls the intensity of the random fluctuations introduced by the Wiener process $W(t)$. The solution to Equation~\eqref{SDE-OU} is given by
\begin{equation}
\label{eq:OU_solution}
\epsilon(t) = \epsilon_0 e^{-\phi t}+\sigma_\text{OU}\int_0^t e^{-\phi(t-\tau)}\text{d}W(\tau),
\end{equation}
where the second term follows $\mathcal{N}(0,\sigma^2(1-e^{-2\phi t})/(2\phi))$. Hence, $\epsilon(t)$ is conditionally normally distributed with
\begin{align}
    \mathbb{E}[\epsilon(t)]&= \mathbb{E}[\epsilon_0]e^{-\phi t}=0,
    \label{OU-expectation}\\
    \mathbb{V}[\epsilon(t)]&=e^{-2\phi t}\mathbb{V}[\epsilon_0 ]+\mathbb{V}\left[\sigma_\text{OU}\int_0^t e^{-\phi(t-\tau)}\text{d}W(\tau)\right]=\frac{\sigma^2_\text{OU}}{2\phi}.
    \label{OU-variance}
\end{align}
Note that the choice of $\epsilon(0)\sim\mathcal{N}(0,\sigma^2_\text{OU} /(2\phi))$ means that $\epsilon(t)$ follows the stationary distribution at all times. 

We use the OU process described above to model autocorrelated observation noise with data collected at observation time points $t_1,\ldots,t_{n_s}$. We set $\bm\epsilon = [\epsilon(t_1), \epsilon(t_2), \ldots, \epsilon(t_{n_s})]$ and use the fact that $\epsilon(t_1)\sim \mathcal{N}(0,\sigma^2_\text{OU} /(2\phi))$, with Equation~\eqref{SDE-OU} describing the noise at subsequent times $t_2,\ldots,t_{n_s}$. Given Equation~\eqref{eq:OU_solution}, it follows that for $i=1,\ldots, n_s-1$ the conditional distribution of $\epsilon(t_{i+1})$ given $\epsilon(t_{i})$ is
\begin{equation}
    \epsilon(t_{i+1}) \mid \epsilon(t_{i}) \sim \mathcal{N}\left(\epsilon(t_i) e^{-\phi \Delta_i}, \frac{\sigma^2_\text{OU}}{2\phi}\left(1 - e^{-2\phi \Delta_i}\right)\right),
\end{equation}
where  $\Delta_i=t_{i+1}-t_i$ is the time between consecutive observations. The autocorrelation is given by $\text{AC}(\Delta_i)=e^{-\phi\Delta_i}$ for time lag $\Delta_i$. Hence we see that the autocorrelation diminishes as the time interval between observations increases, and that the smaller the value of $\phi$, the stronger the autocorrelation. 

The likelihood function can be expressed as
\begin{align}
    \nonumber
    L(\bm \epsilon) & =   -\frac{n_s}{2} \log\left(\pi\frac{\sigma^2_\text{OU}}{\phi}\right) - \sum_{i=1}^{n_s-1} \frac{1}{2} \log\left(  1 - e^{-2\phi \Delta_i} \right) \\
    & \qquad \qquad  -\frac{\phi}{\sigma^2_\text{OU}}\epsilon^2(t_1)-\frac{\phi}{\sigma^2_\text{OU}}\sum_{i=1}^{n_s-1} \frac{\left(\epsilon(t_{i+1}) - \epsilon(t_{i}) e^{-\phi \Delta_i}\right)^2}{ 1 - e^{-2\phi \Delta_i}},
\end{align}
or, equivalently,
\begin{equation}
    L(\bm \epsilon)=-\frac{n_s}{2}\log(2\pi)-
    \frac{1}{2}\log(\mathrm{det}(\bm{\Sigma}))-\frac{1}{2} \bm \epsilon  \bm{\Sigma}^{-1}\bm \epsilon^{\mathrm{T}}, 
\end{equation}
where $\bm{\Sigma} \in \mathbb{R}^{n_s\times n_s}$ is the covariance matrix of the OU process, with entries
\begin{equation}\label{Covariance_OU}
\Sigma_{ij} = 
        \mathrm{Cov}[\epsilon(t_i),\epsilon(t_j)]=\frac{\sigma^2_\text{OU}}{ 2\phi}e^{-\phi |t_i-t_j|},
\end{equation}
for $i,j=1,\ldots,n_s$. Given that $\bm \epsilon=\bm Y-\bm C(\bm \theta)$, we can write
\begin{equation}
\label{log_likelihood_OU_matrix}
    L(\bm{\theta}|\bm{Y})=-\frac{n_s}{2}\log(2\pi)-
    \frac{1}{2}\log(\mathrm{det}(\bm{\Sigma}))-\frac{1}{2}(\bm{Y}-\bm{C}(\bm{\theta}))\bm{\Sigma}^{-1}(\bm{Y}-\bm{C}(\bm{\theta}))^{\mathrm{T}}. 
\end{equation}

To avoid issues with practical identifiability, we assume that $\phi$ is known, and we apply maximum likelihood estimation to estimate the volatility coefficient as
\begin{equation}
\label{eq_sigma_estimate_OU}
    \hat\sigma^2_\text{OU} =\frac{ 2\phi}{n_s}(\bm{Y}-\bm{C}(\bm \theta)) \tilde {\bm\Sigma}^{-1}(\bm{Y}-\bm{C}(\bm \theta))^{\mathrm T},
\end{equation}
where the entries of $\tilde{\bm \Sigma}$ are  
\begin{equation}
\label{tilde_Sigma}
\tilde \Sigma_{ij} = 
         e^{-\phi |t_i-t_j|}.
\end{equation}
We note that $\sigma^2_\text{OU}$ does not affect the process of optimising the experimental design as it serves merely as a scaling factor in the objective function. 

The entries of the Fisher information matrix, $\mathcal{F}$, can be written  
\begin{equation}
\label{FIM_OU_3}
    F_{i,j} =\theta_i\theta_j\frac{\partial \bm{C}}{\partial \theta_i}  \bm\Sigma^{-1}\left(\frac{\partial \bm{C}}{\partial \theta_j}\right)^\mathrm{T},
\end{equation}
and entries of the global sensitivity matrix, $\mathcal{G}$, are defined as
\begin{equation}
\label{def_OUSob}     G_{ij}=\bm{S}_i\bm{\Sigma}^{-1}\bm{S}_j^{\mathrm{T}},
\end{equation}
where $\bm\Sigma$ is as defined in Equation~\eqref{Covariance_OU} and, as before, $\bm{S}_i=[S_i(t_1),\ldots,S_i(t_{n_s})]$ is a row vector composed of the total effect Sobol' indices of parameter $i$ at each time $t_s$, $s=1,\ldots,n_s$, as defined in Equation~\eqref{Sobol_total}. It is worth noting that the inverse of $\bm\Sigma$, also known as the precision matrix, contains information about the conditional dependencies between different time points, and so their inclusion in the sensitivity matrices $\mathcal{F}$ and $\mathcal{G}$ enable quantification of how parameter sensitivities are influenced by temporal correlations in the data. 


\subsection{Optimal experimental design} 
\label{section_OED}

The objective of optimal experimental design in this work is to improve the reliability of parameter estimation through optimising the experimental conditions; in this work, we will focus on optimising the observation times $ t_1,\ldots,t_{n_s}$. 


\paragraph{Fisher information matrix.} In the literature, there are three main quantities calculated from the Fisher information matrix that are routinely used to measure the quality of parameter estimates, namely the trace, the determinant and the minimum expected information gain~\cite{banks2011comparison,bauer2000numerical}. We found these metrics to yield similar results when used for optimizing the selection of observation times. Therefore, we choose to present results generated using the determinant of the Fisher information matrix. 

We denote the set of measurement times as $\bm{t}=[t_1,\ldots,t_{n_s}]$ so that the aim is to optimise
\begin{equation}\label{OED_FIM}
\max_{\bm{t}}~f=\max_{\bm{t}}\left[\vphantom{\sum}\text{det}(\mathcal{F}(\bm{\theta},\bm{t}))\right],
\end{equation}
where $\mathcal{F}(\bm{\theta},\bm{t})$ is the Fisher information matrix evaluated for parameters $\bm\theta$ at times $\bm{t}$. We impose the following constraints on $\bm{t}$:
\begin{equation}
t_1\geq0, \quad t_{n_s}\leq{t_\text{final}}, \quad t_{s+1}-t_s\geq\underline{\Delta}_t \text{\, for \,} s=1,\ldots,n_s-1,
\end{equation}
where $\underline{\Delta}_t$ is the minimum time interval between observations that is practical for experimental design (throughout this work we take $\underline{\Delta}_t=2.0$). This optimal experimental design problem is a non-convex, continuous nonlinear programming problem~\cite{danilova2022recent}. We solve it using the fmincon function in Matlab with the interior-point method. To mitigate the risk of the algorithm converging to local optima,  we implement a random restart technique, which involves generating $50$ different random initial guesses to serve as starting points for fmincon. 


\paragraph{Global information matrix.} Since Sobol' indices need to be calculated using a Monte Carlo method, it is prohibitively time-consuming to consider the same optimisation problem as that for the Fisher information matrix. To make progress, we instead establish a mixed-integer nonlinear programming problem. We first, as for the Fisher information matrix, define $\underline{\Delta}_t$ to be the minimum time interval between observations that is practical for experimental design (again, taking $\underline{\Delta}_t=2.0$ throughout this work). We then define a fine grid of $k$ potential observation times $[t_1^p,\ldots,t_k^p]$, with $t_i^p=\underline{t}+(i-1)\underline{\Delta}_t$ for $i=1,\ldots,k-1$, and $t^p_k=t_\text{final}$, where $k$ is the largest integer such that $t_{k-1}^p\leq{t_\text{final}}$. The optimal experimental design aim then is to select the $n_s<k$ time points that give rise to the optimal objective function, which is now $\text{det}(\mathcal{G})$, where $\mathcal {G}$ is the global information matrix defined in Equation~\eqref{def_IIDSob} (uncorrelated noise) or Equation~\eqref{def_OUSob} (correlated noise).

We formulate the optimisation problem using the decision variables $\bm{I}=[I_1,\ldots,I_k]$, where $I_i\in\{0,1\}$ for $i=1,\ldots,k$ indicates whether an observation is made at time point $t_i$ or not. The aim is then to optimise
\begin{equation}
\label{OED_Sob}
    \max_{\bm{I}}~g=\max_{\bm{I}}\left[\vphantom{\sum}\text{det}(\mathcal{G}(\bm{\theta},\bm{I}))\right],
\end{equation}
subject to the constraint 
\begin{equation}
\label{eq:constraint_1}
    \sum_{i=1}^{k} I_i=n_s.
\end{equation} 
We solve the optimisation problem in Equations~\eqref{OED_Sob}--\eqref{eq:constraint_1} using PlatEMO, an evolutionary multi-objective optimisation platform for Matlab~\cite{tian2017platemo,tian2023practical} and, in particular, we use the competitive swarm optimizer. Throughout this work we use the parameter ranges $r\in[0.14,0.26]$, $K\in[35,65]$ and $C_0\in[3.15,5.85]$.


\subsection{Profile likelihood approach} 
\label{section_Profile}

In order to assess the performance of the observation time points selected by solving the optimal experimental design problem, we calculate the profile likelihood. The practical identifiability of a parameter can be determined by examining the shape and width of the profile likelihood function~\cite{raue2009structural}. A parameter is considered practically identifiable if its profile likelihood is well-defined, meaning it is unimodal with a single peak, and a finite confidence interval. Conversely, a flat or broad profile likelihood suggests that a wide range of parameter values yield similar fits to the data, thereby indicating that the parameter is not practically identifiable. 

To calculate the univariate profile likelihoods, we partition the parameter vector $\bm{\theta}$ into the interest parameter, $\theta_i$, and nuisance parameters, $\bm{\theta}_{\sim{i}}$. We define the univariate, normalised profile log-likelihood function for the interest parameter $\theta_i$,  
\begin{equation}
\label{profile-like}    
l_p(\theta_i|\bm{Y})=\sup_{\bm{\theta}_{\sim i}}L(\bm{\theta}|\bm{Y})-\sup_{\bm{\theta}}L(\bm{\theta}|\bm{Y}).
\end{equation}
Given this definition, we can use $l_p(\theta_i|\bm{Y})=1.92$ to define an approximate 95\% confidence interval for $\theta_i$~\cite{royston2007profile}, as follows
\begin{equation}\label{eq:definition_CI}
    \mathrm{CI}=\left\{\theta_i|l_p(\theta_i|\bm{Y})\ge -\chi^2_{0.95;1}/2 = -1.92\right\},
\end{equation}
where $\chi^2_{0.95;1}$ represents the critical value of the chi-square distribution with one degree of freedom at the $95\%$ confidence level. We use the fmincon function in Matlab to calculate the MLE, the profile likelihood and 95\% confidence interval for each parameter.

We use the 95\% confidence intervals to generate prediction intervals for the model~\cite{simpson2024parameter}. The prediction intervals are calculated through parameter sampling and log-likelihood evaluation: parameters are uniformly sampled and, for each sample, the normalised log-likelihood is computed and compared to a threshold of $\chi^2_{0.95;3}/2$. Parameters that exceed this threshold are retained, ensuring they lie within the $95\%$ confidence region. From the valid samples, $M$ population trajectories are generated at time points $t=1,\ldots,n_s$, and the upper and lower bounds of these trajectories are calculated. Noise model bounds, based on the $5\%$ and $95\%$ quantiles of the relevant distribution, are then added into the upper and lower bounds, respectively, so as to form the prediction interval.


\section{Results}
\label{section_results}

In this section we present the optimal experimental design results and the impact of optimally choosing measurement times on the parameter estimates. 


\subsection{Parameter identifiability under the different noise models}
\label{section3_identifiability}


\begin{figure}[htb!]    
    \begin{center}
    \includegraphics[width=\textwidth]{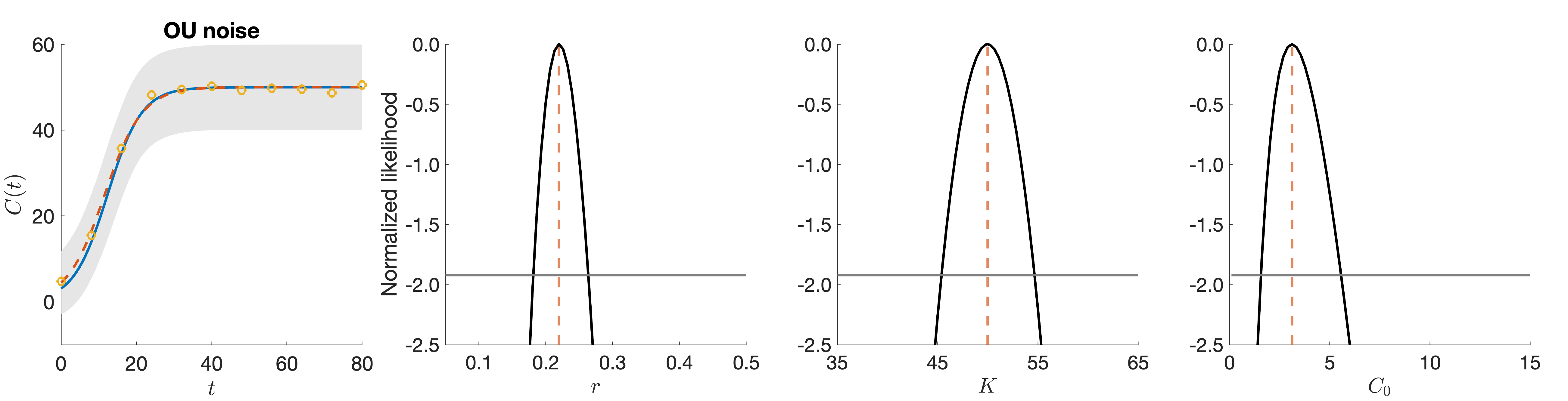}\\
    \includegraphics[width=\textwidth]{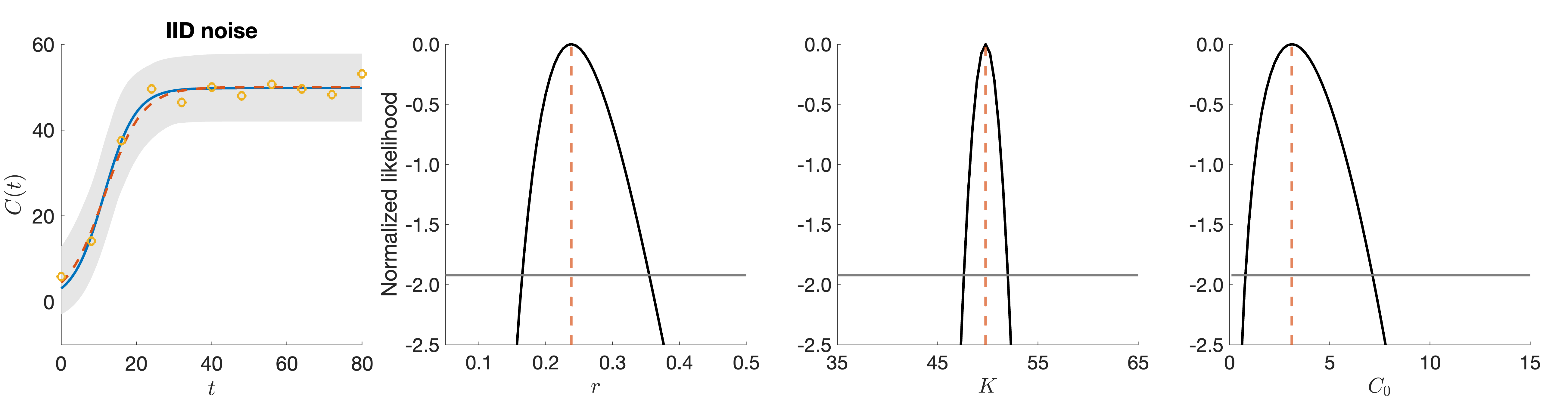}\\
    \includegraphics[width=\textwidth]{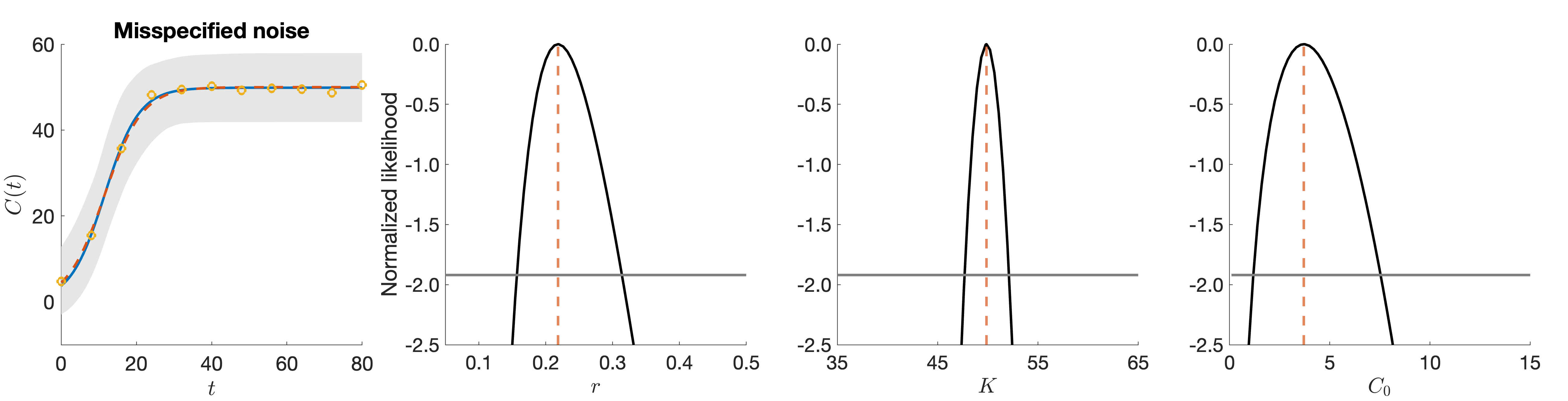}
    \end{center}
	\caption{Parameter identifiability under the different noise models, with data collected at 11 regularly spaced time points $t=0,8,16,\ldots,80$. The left column shows the sampled data, the underlying model solution using the input parameters, and the solution generated using the MLE parameters. The shaded region shows the 95\% prediction interval generated using the 95\% confidence intervals from the profile likelihoods. The right-most three panels show the profile likelihoods for each parameter. The results in the top row are generated using uncorrelated Gaussian noise with $\sigma_{\rm{IID}}^2=9.0$. The results in the middle row are generated using correlated OU noise with $\phi=0.02$ and $\sigma_{\rm{OU}}^2/(2\phi)=9.0$. The results in the bottom row are generated using a misspecified noise model: the data are generated using correlated OU noise with $\phi=0.02$ and $\sigma_{\rm{OU}}^2/(2\phi)=9.0$, whilst the analysis is carried out assuming IID Gaussian noise with $\sigma_{\rm{IID}}^2=9.0$.}
	\label{fig:posterior_profile}
\end{figure}


We first explore the extent to which the model parameters can be estimated under different noise models when the data are sampled at regular time intervals, with 11 observations at times $t=0,8,16,\ldots,80$. Results for $\phi=0.02$ are shown in Figure~\ref{fig:posterior_profile}, with similar results for $\phi=0.1$ (less correlated noise) shown in Supplementary Figure S1. The left column shows the sampled data (circles) along with the underlying model solution using the input parameters (dashed line) and the solution generated using the MLE parameters (solid line), and the shaded region shows the 95\% prediction interval. The top row shows the results from the uncorrelated noise model, the second row shows the results from the correlated noise model, whilst the bottom row shows the results in the case of misspecification in the noise model, where correlated noise was used to generate the observations but the profile likelihoods were generated using the uncorrelated noise model. To enable a comparison of the effects of uncorrelated observation noise and correlated observation noise on parameter inference under consistent variance, we choose parameters such that  $\sigma^2_{\rm{IID}}=\sigma_{\rm{OU}}^2/(2\phi)$. We see that in all cases all model parameters are practically identifiable, with the widths of the confidence intervals depending on the both the noise model and the data. A key point to note is that incorrectly assuming uncorrelated noise can impact the parameters in different ways. For example, the carrying capacity, $K$, is predicted to be more confidently estimated, whilst the uncertainty in the growth rate is predicted much higher in the misspecified case. 


\begin{figure}[tb]
    \begin{center}
    \includegraphics[width=0.95\textwidth]{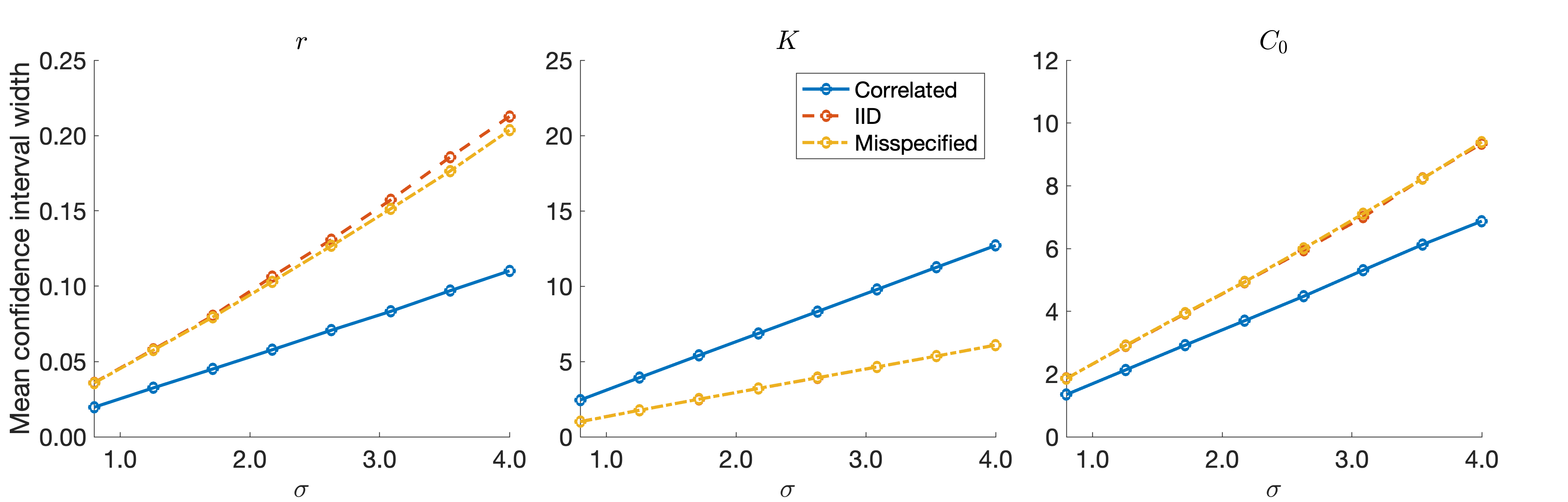}
    \includegraphics[width=0.95\textwidth]{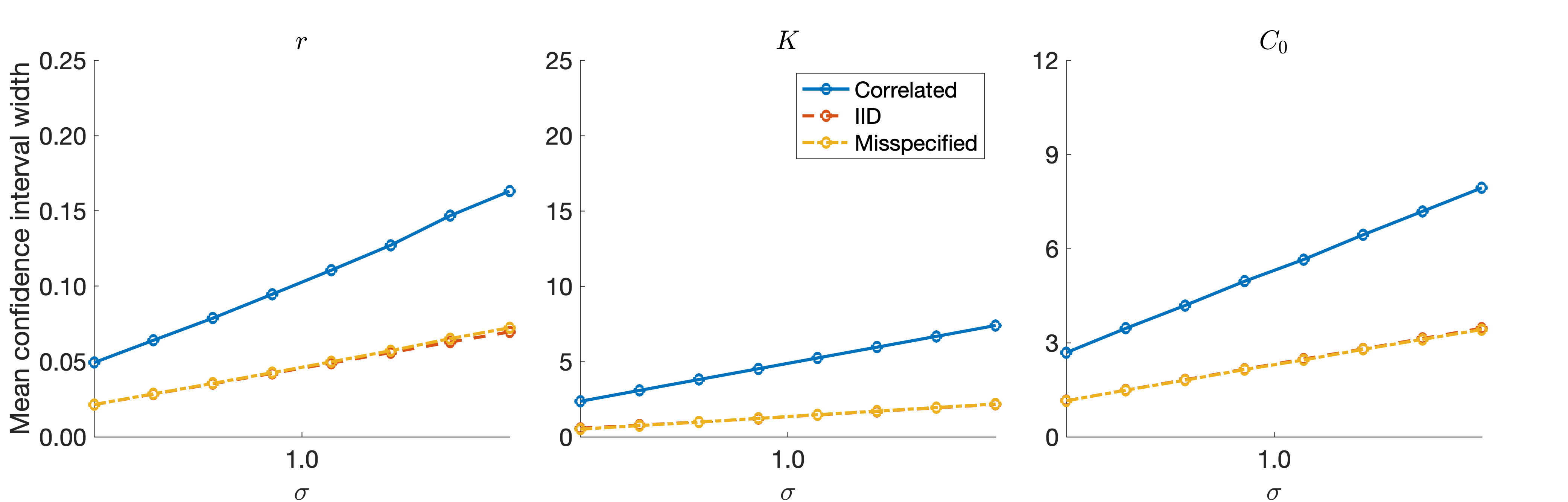}
    \end{center}
    \caption{Change in the mean $95\%$ confidence interval for parameters $r$, $K$ and $C_0$ as the variance increases under both correlated and uncorrelated noise, and when the noise is misspecified (as per Figure~\ref{fig:posterior_profile}). Top row: $\sigma^2=\sigma_\text{IID}^2=\sigma_\text{OU}^2/(2\phi)$ with $\phi=0.02$ so that the variance is equivalent across the different noise models. Bottom row: $\sigma^2=\sigma_{\rm{IID}}^2=\sigma_{\rm{OU}}^2$ with $\phi=0.02$ so that the variance is larger in the correlated noise process than the uncorrelated noise process. In all cases, results were generated by averaging the confidence interval width from 1000 simulations for each noise model and parameter set, and 11 observations were used.}
    \label{fig:CI_varysigmaC}
\end{figure} 


Figure~\ref{fig:CI_varysigmaC} shows how the width of the confidence intervals changes with the variance for both uncorrelated and correlated noise, as well as misspecified noise\footnote{Recall that in the misspecified noise case, the data are generated using correlated data but the confidence intervals are generated assuming the data are uncorrelated.}, in each case using 11 regularly spaced observation time points, $t=0,8,16,\ldots,80$. The top row shows results with $\sigma^2=\sigma_\text{IID}^2=\sigma_\text{OU}^2/(2\phi)$, and the mean reversion rate of the OU process held constant at $\phi=0.02$ so that the noise variance is equivalent across both noise models. We see that, for all values of $\sigma^2$, the confidence intervals for all parameters are almost identical for uncorrelated noise and misspecified noise, with greater confidence in the parameter estimated for the carrying capacity $K$, and less confidence in the estimates for the growth rate, $r$, and the initial condition, $C_0$, compared to the correlated noise case. This plots highlights the potential pitfalls of incorrect assumptions on the noise model, in that it can potentially lead a modeller to either have too little, or too much, confidence in parameter estimates. The bottom row of Figure~\ref{fig:CI_varysigmaC} illustrates that setting $\sigma_{\rm{IID}}^2=\sigma_{\rm{OU}}^2$ with significant noise correlation ($\phi=0.02$) results in greatly reduced confidence in the parameter estimates for the correlated noise case compared to the uncorrelated and misspecified cases. This is to be expected as, in this case, the variance of the correlated noise process is 25 times larger than that of the uncorrelated noise process. 


\subsection{Optimal experimental design}
\label{optimisation_Results}

In this section, we use the framework outlined in Section~\ref{section_OED} to explore how the optimal time points for observations change as the total number of observations, $n_s$, is varied, and also how correlations in the observation noise influence the optimal measurement protocol. We use both the local measures obtained using the Fisher information matrix, and global measures obtained using the total Sobol' index, as outlined in Section~\ref{section_information}. All results are generated using $\phi=0.02$, which entails a significant degree of correlation in the noise process. 


\begin{figure}[tb]    
    \begin{center}
    \includegraphics[width=\textwidth]{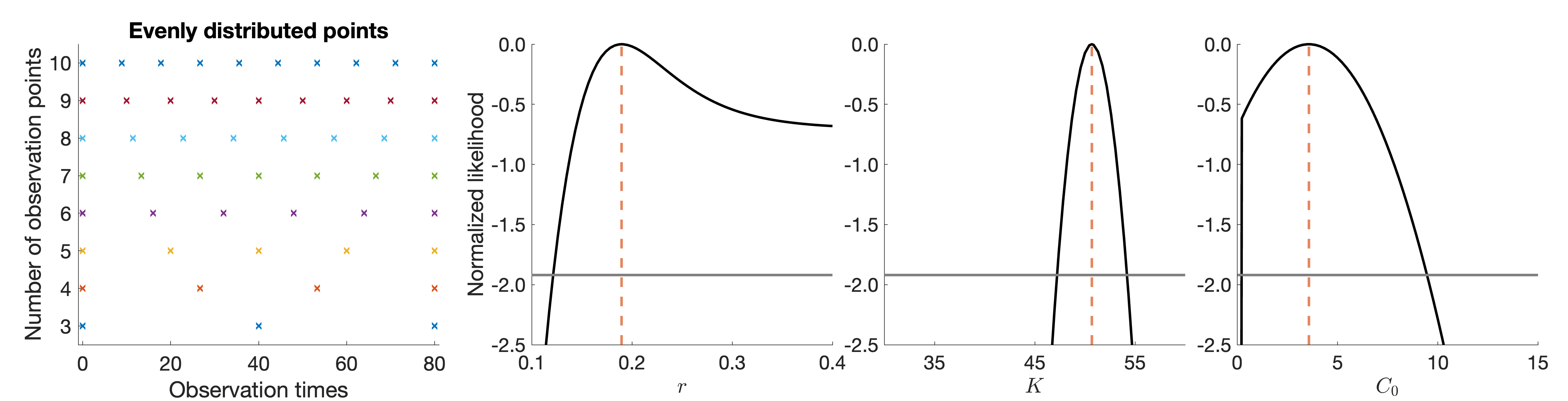}\\
    \includegraphics[width=\textwidth]{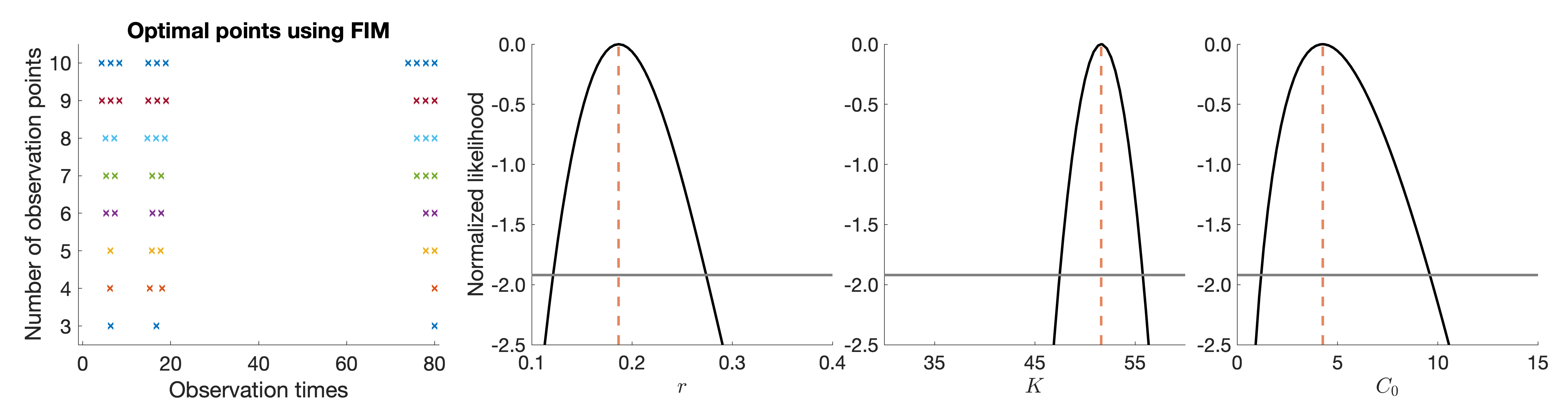}\\
    \includegraphics[width=\textwidth]{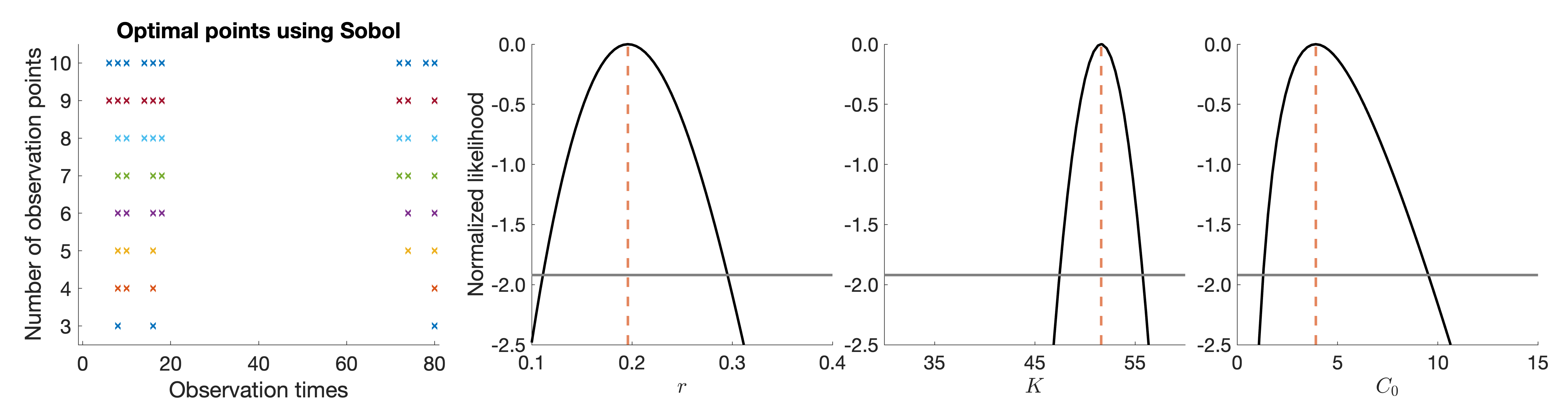}\\
    \end{center}
	\caption{The results of optimising experimental design in the uncorrelated noise case as the number of observations, $n_s$, is varied from three to ten. The top row shows the results from evenly distributed observations, whilst the second row shows the optimal experimental design derived using the Fisher information matrix, and the bottom row shows the optimal experimental design derived using the global information matrix. In each case, the left-hand plots show the observations time points, whilst the right-hand three plots show the profile likelihoods obtained from five observations. In all plots, $\sigma_{\rm{IID}}^2 = 9.0$.}	
    \label{fig:opt_points_IID}
\end{figure}

\begin{figure}[tb]
    \hspace{2cm}
    \includegraphics[width=\textwidth]{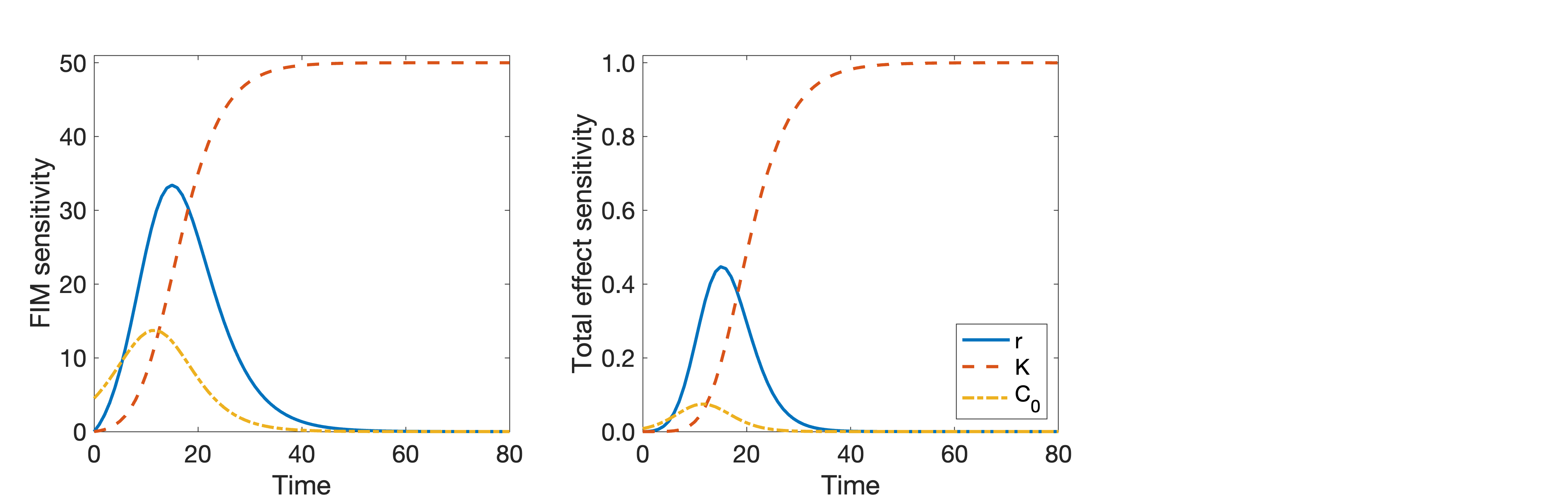}
    \caption{(a) Fisher information sensitivity and (b) global information sensitivity, formulated in terms of the gradient of the output, $C(t;\bm{\theta})$, with respect to the individual parameters. We evaluate the sensitivity at $t_k=0,2,4,\ldots,80$, with the Fisher information sensitivity calculated as $\theta_i\,\partial C(t_k)/\partial\theta_i$ and the global information sensitivity calculated as $S_i(t_k)$, as defined in Equation~\eqref{Sobol_total}.}
    \label{fig:Sensitivity_index}
\end{figure}

Figure~\ref{fig:opt_points_IID} shows the results of optimising the experimental design for uncorrelated noise, detailing the placement of the observations as the total number of observations is varied from three to ten. It also compares the confidence in parameter estimates under both optimised and evenly spread observations when five observations are made in total. We see that both the Fisher information matrix and the global information matrix lead to experimental designs that place the observations into three groups, with the first two groups of observations lying in the first 20 time units, which corresponds to the exponential growth phase, and the third group gathered towards the final time, which corresponds to the process being at steady state (carrying capacity). These results make sense intuitively, with the first two groups of time points enabling accurate inference of the growth rate, $r$, and the initial condition, $C_0$, and the third group enabling accurate inference of the carrying capacity, $K$. These intuitive explanations are further supported by the plots in Figure~\ref{fig:Sensitivity_index} which shows how the sensitivity of the Fisher information matrix and global information matrix measures to variation in each parameter change over the course of the experiment. Our results also highlight how careful placement of the observations can enable accurate inferences to be made with fewer time points: the right-hand three plots highlight that under the optimised protocols, all parameters can be accurately estimated using just five observations, whereas neither the growth rate, $r$, or the initial condition, $C_0$, can be accurately inferred using five equally spaced observations. 


\begin{figure}[tb]    
	\begin{center}
    \includegraphics[width=\textwidth]{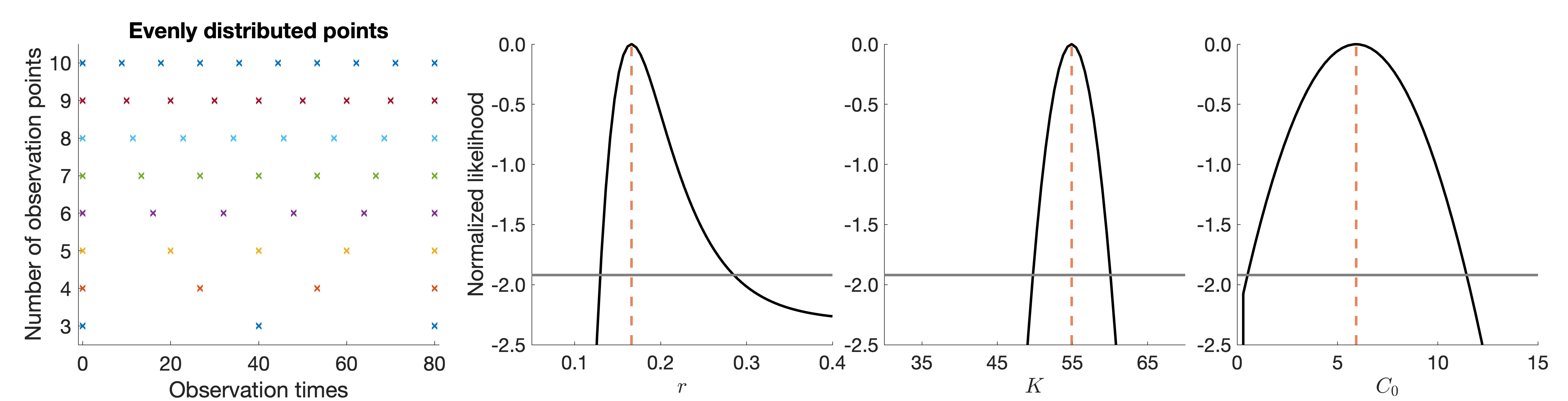}\\
    \includegraphics[width=\textwidth]{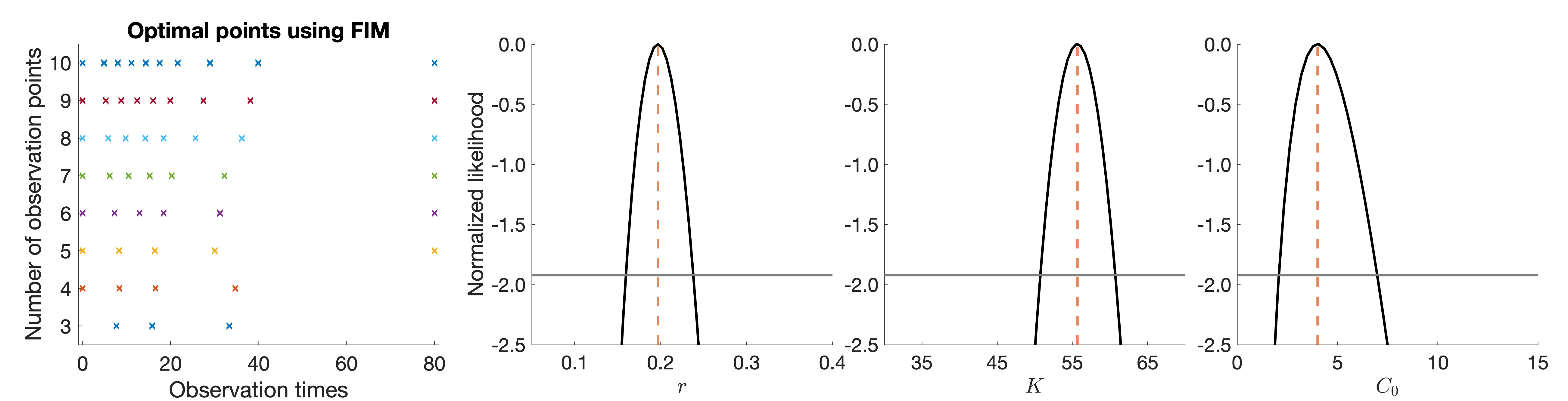}\\
    \includegraphics[width=\textwidth]{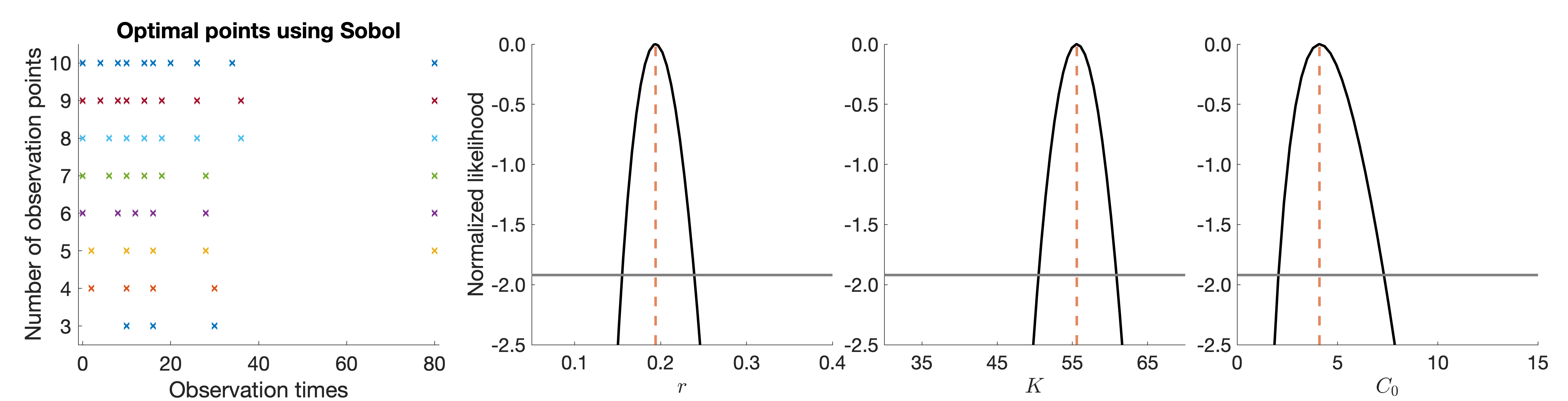}\\
    \end{center}
	\caption{The results of optimising experimental design in the correlated noise case as the number of observations, $n_s$, is varied from three to ten. The top row shows the results from evenly distributed observations, whilst the second row shows the optimal experimental design derived using the Fisher information matrix, and the bottom row shows the optimal experimental design derived using the global information matrix. In each case, the left-hand plots show the observations time points, whilst the right-hand three plots show the profile likelihoods obtained from five observations. In all plots, $\sigma_{\rm{OU}}^2/(2\phi)=9.0$.}
    \label{fig:opt_points_OU}
\end{figure}

When the observation noise is correlated (see Figure~\ref{fig:opt_points_OU}), the optimisation results change dramatically: the majority of the observation points are evenly distributed during the first half of the experiment, until approximately 40 time units, with only a single observation in the second half of the experiment, which is generally placed at, or very close to, the terminal time, $t_\text{final}$. This change in experimental design is intuitive: the significant correlation in the noise process (small $\phi$) means that less information about the parameters is gained from two closely placed observations. Once again, we see that the Fisher information matrix and the global information matrix provide very similar results in terms of optimal placement of the observation time points. We also note that with correlated observation noise, all parameters can in fact be inferred from five observations without any optimisation---evenly spacing the time points provides closed confidence intervals for all model parameters---although optimisation clearly reduces the uncertainty in the parameter estimates. 


\subsection{Impact of optimal experimental design upon parameter identifiability} 
\label{Result_test}

\begin{figure}[tb]    
    \begin{center}
    \includegraphics[width=0.95\textwidth]{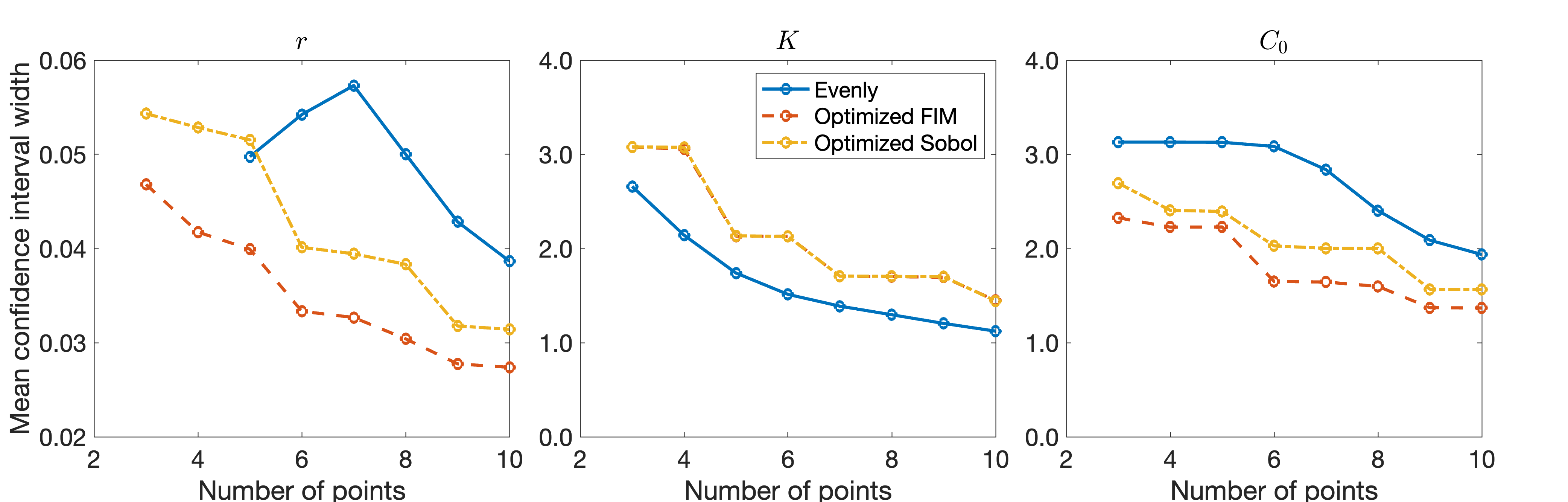}\\
    \includegraphics[width=0.95\textwidth]{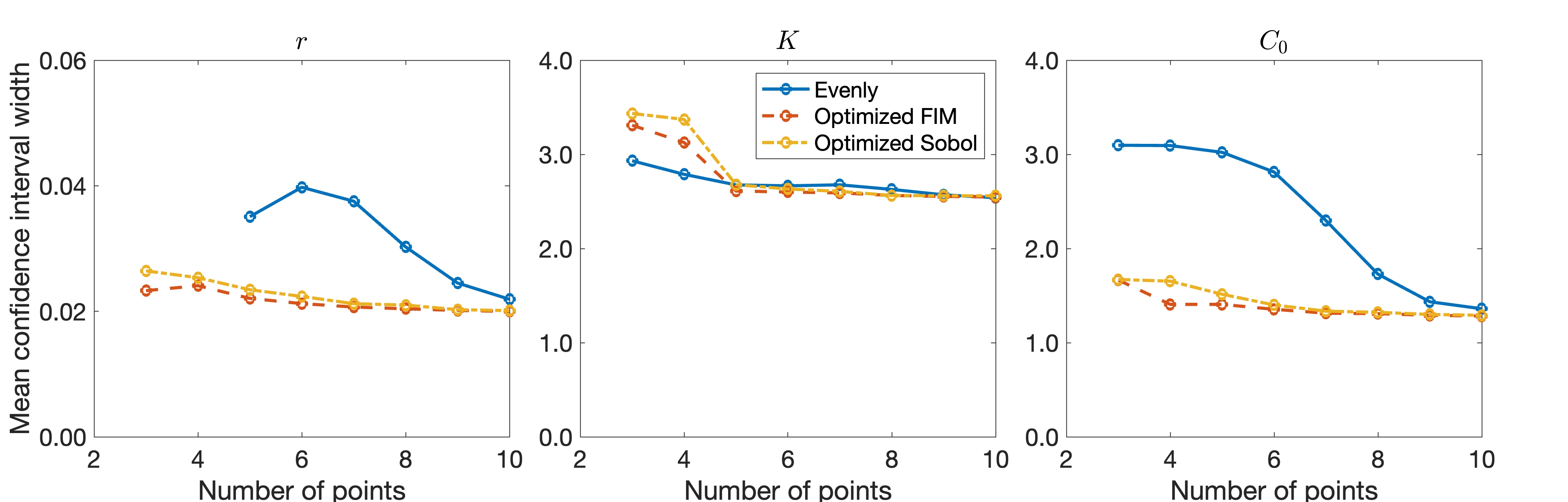}\\
    \includegraphics[width=0.95\textwidth]{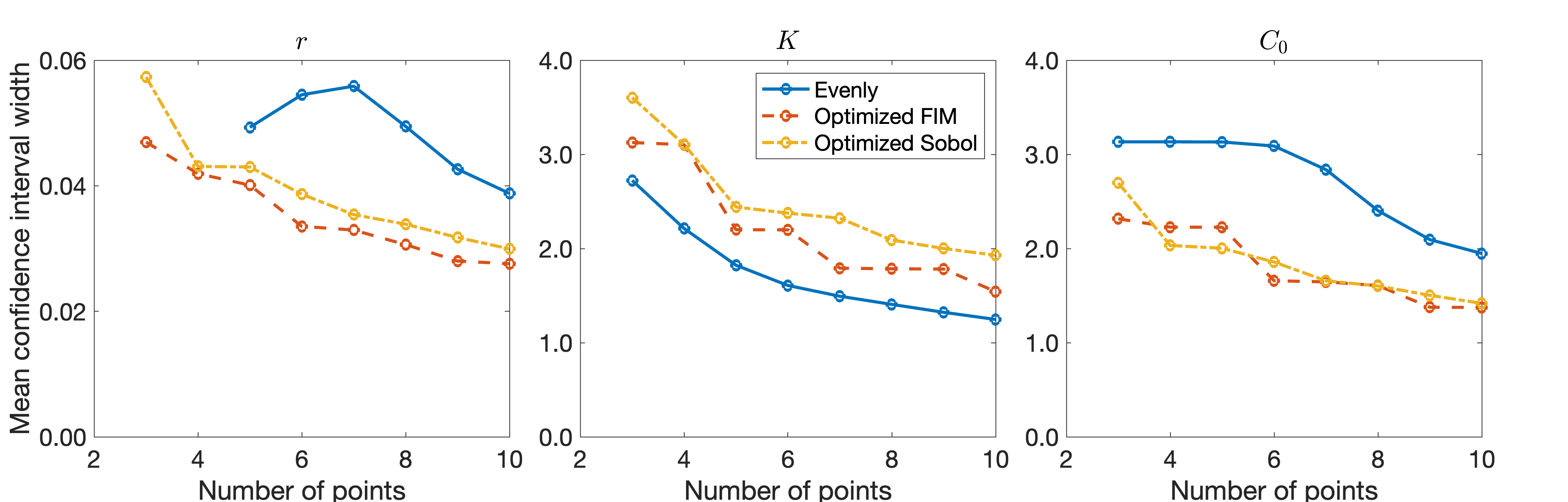}\\
    \end{center}
    \caption{Mean parameter confidence intervals in the case of uncorrelated noise (top row), correlated noise (middle row) and misspecified noise (bottom row). In each case, we plot the mean confidence interval as a function of the number of observation points, $n_s$. The plots in each row are generated used 1000 simulations with $\sigma_{\rm{IID}}=0.8$ (top row), $\sigma_{\rm{OU}}=0.16$ and $\phi=0.02$ (middle and bottom rows) so that $\sigma_\text{IID}^2=\sigma_\text{OU}^2/(2\phi)=0.64$.}
    \label{fig:CI_number_samples}    
\end{figure}

In this section, we explore how changes in the number of observation time points impact the confidence in parameter estimates. The mean parameter confidence intervals, calculated using the profile likelihood method presented in Section~\ref{section_Profile}, are shown in Figure~\ref{fig:CI_number_samples}. The top row displays results for the uncorrelated noise model, the second row for the correlated noise model, and the bottom row for the case of model misspecification, where correlated noise was used to generate observations, but profile likelihoods were computed assuming an uncorrelated noise model. Additional results, for different values of $\phi$, are shown in Supplementary Information Figure S2.


Figure~\ref{fig:CI_number_samples} shows the intuitive result that, in all cases, parameter estimates become more confident as the number of observation time points is increased, and that optimising the timing of observations generally improves the confidence in parameter estimates compared to the case of evenly distributing the observation time points. The exception to this rule is the carrying capacity parameter, $K$, where the confidence intervals are often marginally smaller for evenly spaced time points. This is a direct result of the fact that the optimisation algorithms place the majority of the observations at early times in order to accurately estimate the growth rate, $r$, and the initial condition, $C_0$. Note that for small numbers of observation time points ($n_s=3,4$) the confidence intervals are half open when the observation points are evenly distributed, hence we do not display a result in those cases. In addition,  note that when the number of evenly distributed points increases from five to seven in Figure~\ref{fig:CI_number_samples}, the width of the confidence interval for the parameter $r$ increases. One possible explanation is that among the five evenly distributed points, some provide more information for estimating the parameter $r$. Close investigation reveals that this is indeed the case: observations placed close to $t=20$ are important for the estimation of $r$ (see Supplementary Figure S3). It is worth noting that for the uncorrelated noise model, the mean confidence intervals for the carrying capacity, $K$, generated using local and global approaches to optimising the observation time points overlap. This is because the selected observation points lie in the saturation region and are nearly identical across the methods (as shown in the left column of Figure 3).  


\subsection{Impact of the autocorrelation level on experimental design} 
\label{Result_correlated}


\begin{figure}[tb]    
    \begin{center}
    \includegraphics[width=0.95\textwidth]{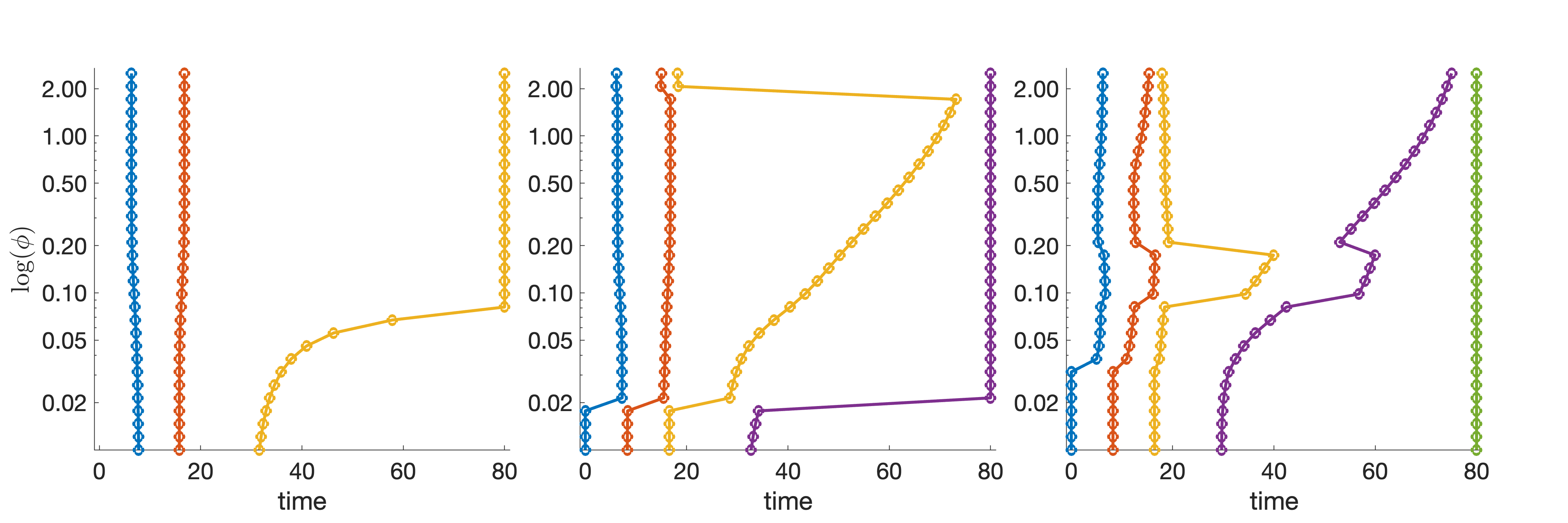}\\
    \includegraphics[width=0.95\textwidth]{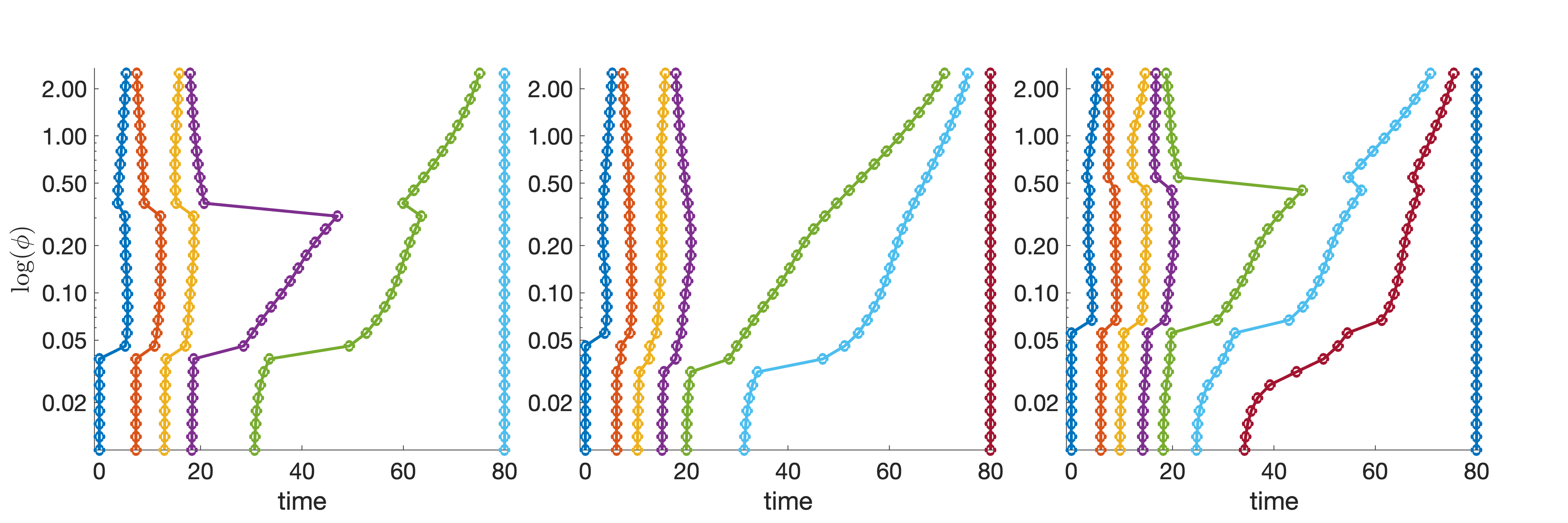} 
    \end{center}
    \caption{The optimal observation time points as a function of the mean-reversion rate, $\phi$. Note that the degree of autocorrelation decreases as $\phi$ increases.}
    \label{fig:opt_phi_vary}
  \end{figure} 


In this section, we will investigate how the mean-reversion rate $\phi$, reflecting the level of autocorrelation in the observation noise, influences the optimised distribution of observation times. Figure~\ref{fig:opt_phi_vary} shows how the optimised observation times (generated using the Fisher information matrix) change as $\phi$ increased (which corresponds to decreasing the autocorrelation). For very high autocorrelation, the majority of measurement points are placed within the first half of the experiment in order that the growth rate, $r$, and initial condition, $C_0$, can be accurately measured. As the degree of autocorrelation diminishes ($\phi$ is increased), a measurement point is placed at, or close to, the terminal time, which allows estimation of the carrying capacity, $K$. With more total observations, $n_s$, further points are added at the end of time interval for larger values of $\phi$, which enables increases in the accuracy of estimation of the carrying capacity, $K$.


\begin{figure}[tb]
    \begin{center}
    \includegraphics[width=0.95\textwidth]{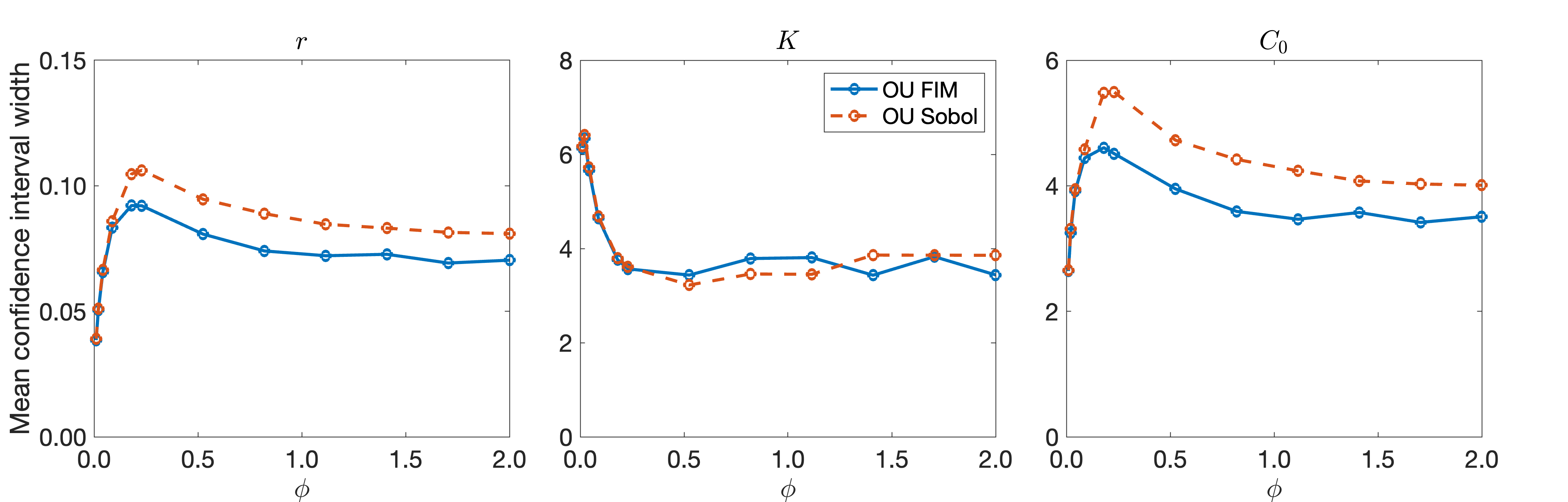}\\
    \includegraphics[width=0.95\textwidth]{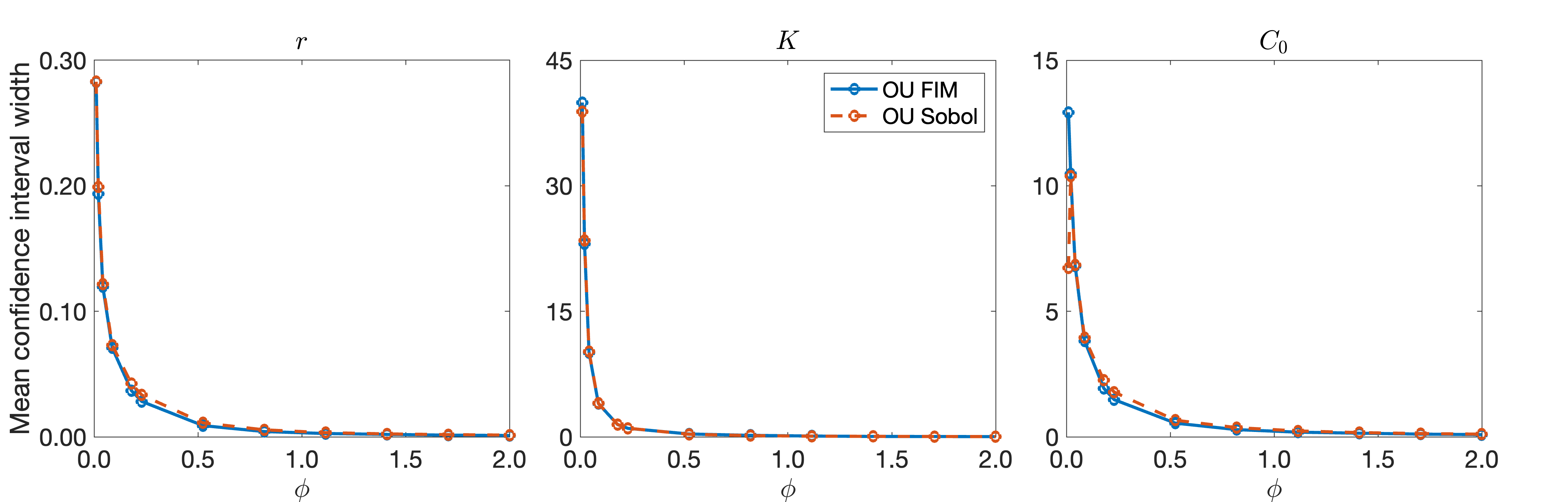}
    \end{center}
    \caption{The mean confidence interval width for parameters estimated using 11 optimised observation time points as the mean-reversion rate, $\phi$, increases from $0.01$ to $2.0$. In the top row the variance is held constant at $\sigma_{\rm{OU}}^2/(2\phi)=4.0$ whilst in the bottow row the variance changes with $\sigma_{\rm{OU}}=0.3$ held constant. In all cases, results were generated using 1000 simulations.}
    \label{fig:CI_phivary}
\end{figure} 


\enlargethispage{0.7cm}

Figure~\ref{fig:CI_phivary} demonstrates how the parameter confidence intervals vary as a function of the degree of autocorrelation. In the top row the variance, $\sigma_{\rm{OU}}^2/(2\phi)$, is held constant as $\phi$ is varied and we see that there is a complicated relationship between the confidence interval width and $\phi$. For very low values of $\phi$, which correspond to highly correlated observation noise, the majority of the observation points are placed early in the experiment which leads to accurate estimates for the growth rate, $r$, and carrying capacity, $C_0$, and far less confidence in estimates of the carrying capacity, $K$. As the autocorrelation decreases (with increasing $\phi$), increased numbers of observations are placed in the second half of the experiment and confidence in estimates of the carrying capacity, $K$, increase significantly. The bottom row demonstrates the intuitive result that, for $\sigma_{\rm{OU}}=0.3$ held constant and increasing $\phi$ (which corresponds to decreasing the variance, $\sigma_{\rm{OU}}^2/(2\phi)$, of the noise process), the mean confidence interval widths for all parameters decreases.


\section{Discussion} 
\label{section_discussion}

This paper investigates the influence of different types of observation noise on optimal experimental design, where the aim is to determine the optimal observation times for accurate and confident parameter inference. Specifically, we analyze both uncorrelated (IID, Gaussian) and correlated (OU) noise and apply local techniques based on the Fisher information matrix, as well as global sensitivity approaches based on Sobol' indices, to optimise observation times for the the logistic growth model. 

Our results highlight several key points. Firstly, the uncertainty in parameter estimates depends on the type of observation noise, with some model parameters becoming more identifiable, and others less so, as the noise model is varied. Second, we demonstrate that, under both uncorrelated and correlated noise, optimising the placement of observation time points reduces the uncertainty in parameter estimates compared to na\"ive placement of measurement points and that, in many cases, this allows the practitioner to confidently estimate parameter values using fewer observations. Thirdly, our results show that  the optimal observation time points are sensitive to the degree of noise correlation, with high autocorrelation favouring the majority of time points in the first half of the experiment. 

On the whole, are results are relatively insensitive to the choice of objective function, though we highlight that use of the Fisher information matrix requires selection of parameter values \textit{a priori} where as the global information matrix takes parameter ranges as inputs, which means that it may be more appropriate for cases in which parameter values are relatively uncertain. The trade-off is the increased computational complexity of the global method compared to that of the local method.

Our methodology is very general, in the sense that it can be applied in any context where it is possible to explicitly write down the likelihood and evaluate the associated Fisher and global information matrices. For example, it can be easily applied in the context of ordinary and partial differential equation-based models, and with a huge range of observation noise models. It would also be possible to integrate a cost-benefit analysis into the framework, for example using multi-objective optimisation to balance the costs of each observation against the quality of parameter estimates. 

In the future, we aim to extend our methodology to other models in mathematical biology to uncover more general insights regarding the effects of autocorrelated measurement processes on parameter estimates for differential equation-based models. Developing a method to diagnose whether noise is correlated or independent would be beneficial. Additionally, it would be promising to optimize other experimental conditions, such as external inputs, beyond just measurement points, or to consider multi-objective optimal experimental design problems.


\section*{Acknowledgements}

J.Q. would like to thank the Mathematical Institute, University of Oxford, for their support and hospitality during a visit to Oxford. R.E.B. acknowledges support from the Simons Foundation (MP-SIP-00001828).


\bibliographystyle{RS}
\bibliography{00-main}


\end{document}